\documentclass[12pt]{article}

\newtheorem{theorem}{Theorem}[section]
\newtheorem{corollary}[theorem]{Corollary}

\newtheorem{condition}[theorem]{Condition}
\newtheorem{lemma}[theorem]{Lemma}
\newtheorem{example}[theorem]{Example}
\newtheorem{proposition}[theorem]{Proposition}
\newtheorem{remark}[theorem]{Remark}

\def\qed{\hbox{}\nobreak [\kern-.4mm] \par \goodbreak \smallskip}

\title{{\bf Logarithm-free $A$-hypergeometric 
 series}}
\author{Mutsumi Saito}

\date{December 28, 2000}

\begin{document}

\maketitle

\begin{abstract}
We give a dimension formula
for the space of logarithm-free series solutions
to an {\it $A$-hypergeometric} 
(or a {\it GKZ hypergeometric}) {\it system}.
In the case where the convex hull spanned by $A$ is a simplex,
we give a rank formula for the system, 
characterize the exceptional set, and
prove the equivalence of
the Cohen-Macaulayness of the toric variety defined by $A$
with the emptiness of the exceptional set.
Furthermore we classify
$A$-hypergeometric systems
as analytic ${\cal D}$-modules.


\end{abstract}

\section{Introduction}

Given a finite set $A$ of integral vectors 
on a hyperplane off the origin and given a parameter vector,
Gel'fand, Kapranov and Zelevinskii \cite{GZK-def} 
defined a system of differential
equations,
called an 
{\it $A$-hypergeometric} (or a {\it GKZ hypergeometric}) {\it system}.
The rank of an $A$-hypergeometric system is greater than or
equal to the volume of the convex hull ${\rm conv}(A)$ spanned by $A$
(Theorem 3.5.1 in \cite{sst-book}).
The set of parameters where the rank is strictly greater than the volume
is called the {\it exceptional set}.
In this paper, we give a dimension formula
for the space of logarithm-free series solutions
to an $A$-hypergeometric system.
Furthermore, in the case where ${\rm conv}(A)$ 
is a simplex,
we give a rank formula, characterize the exceptional set, and
prove the equivalence of
the Cohen-Macaulayness of the toric variety defined by $A$
with the emptiness of the exceptional set.

In the paper \cite{IsoClass}, we defined a finite set $E_\tau(\beta)$
associated to a parameter $\beta$ and a face $\tau$ of the cone 
generated by $A$,
to classify parameters according to isomorphism classes of
their corresponding algebraic $A$-hypergeometric systems.
In this paper, we define a finite set $E_\tau(\beta)$ 
in the same way as in \cite{IsoClass},
associated to a parameter $\beta$ and
a face $\tau$ of the regular triangulation $\Delta_w$ of
${\rm conv}(A)$
determined by a generic weight vector $w$,
to investigate logarithm-free $A$-hypergeometric series.

Roughly speaking, each element of the set $E_\tau(\beta)$ gives
volume-of-$\tau$-many linearly independent logarithm-free
series solutions converging in the direction of $w$ (Theorems 
\ref{theorem:GeneralConvergence} and
\ref{theorem:convergence}).
Using this, we prove that if $E_\tau(\beta)\not= E_\tau(\beta')$ for
some face $\tau$, then 
the analytic $A$-hypergeometric systems with parameter $\beta$ and with
parameter $\beta'$ are not isomorphic (Theorem \ref{thm:AnalyticVersion}).
By counting canonical series solutions associated with $\tau$
but not with smaller faces, we obtain a dimension formula
for the space of logarithm-free
series solutions converging in the direction of $w$ 
(Theorem \ref{theorem:LogFreeFormula}).

When the convex hull ${\rm conv}(A)$ is a simplex, we take a generic
weight vector $w$ such that $\Delta_w$ coincides with ${\rm conv}(A)$
itself. Then the faces of the cone generated by $A$ and those of
$\Delta_w$ have a natural one-to-one correspondence.
In this situation, we prove that all canonical series solutions with
respect to $w$ are logarithm-free 
(Propositions \ref{proposition:FakeIndicial} and \ref{proposition:Exponents}).
Hence the formula in Theorem \ref{theorem:LogFreeFormula} 
is also a rank formula
for $A$-hypergeometric systems (Theorem \ref{simplex:log-free-sol}).
From the rank formula we characterize the exceptional set 
(Theorem \ref{theorem:ExceptionalSet}),
and then we prove the equivalence of the Cohen-Macaulayness of the affine
toric variety defined by $A$ with the emptiness of the exceptional set
(Theorem \ref{theorem:CMequiv}).

Matusevich \cite{Laura} has recently proved the equivalence 
in the codimension 2 case,
and has given examples 
in which the exceptional sets are infinite.

\section{Canonical $A$-hypergeometric series}

In this section, we recall logarithm-free canonical $A$-hypergeometric
series. For details, see \cite{sst-book}.

Let $A=(a_1,\ldots, a_n)=(a_{ij})$ be a $d\times n$-matrix
of rank $d$ with coefficients in ${\bf Z}$. 
Throughout this paper,
we assume that all $a_j$ belong to one hyperplane off the origin 
in ${\bf Q}^d$.
Let ${\bf N}$ be the set of nonnegative integers, and 
${\bf k}$ a field of characteristic zero.
Let $I_A$ denote the toric ideal in the polynomial ring
${\bf k}[\partial]={\bf k}[\partial_1,\ldots,\partial_n]$, i.e.,
$$
I_A=\langle \partial^u-\partial^v\, :\, Au=Av, \, u, v\in {\bf N}^n
\rangle \subset {\bf k}[\partial].
$$
Here and hereafter we use the multi-index notation;
for example, $\partial^u$ means $\partial^{u_1}\cdots \partial^{u_n}$
for $u={}^t(u_1,\ldots,u_n)$.
Given a column vector $\beta={}^t (\beta_1,\ldots,\beta_d)\in {\bf k}^d$,
let $H_A(\beta)$ denote the left ideal of the Weyl algebra
$$
D={\bf k}\langle x_1,\ldots, x_n,\partial_1,\ldots,
\partial_n\rangle
$$
generated by $I_A$ and $\sum_{j=1}^n a_{ij}\theta_{j} -\beta_i$
($i=1,\ldots, d$),
where $\theta_j =x_j\partial_j$.
The quotient $M_A(\beta)=D/H_A(\beta)$ is called
the {\it $A$-hypergeometric system with parameter $\beta$},
and a formal series annihilated by $H_A(\beta)$
an {\it $A$-hypergeometric series with parameter $\beta$}.

Fix a generic weight vector $w\in {\bf R}^n$.
The ideal of the polynomial ring ${\bf k}[\theta]=
{\bf k}[\theta_1,\ldots,\theta_n]$ defined by
\begin{equation}
\label{def:fin}
\widetilde{{\rm fin}}_w(H_A(\beta))
:=D\cdot {\rm in}_w(I_A)\cap {\bf k}[\theta]+\langle A\theta-\beta\rangle
\end{equation}
is called the {\it fake indicial ideal},
where ${\rm in}_w(I_A)$ denotes the initial ideal of $I_A$ with respect
to $w$, and
$\langle A\theta-\beta\rangle$ denotes the ideal generated by
$\sum_{j=1}^n a_{ij}\theta_j -\beta_i$ ($i=1,\ldots, d$).
Each zero of $\widetilde{{\rm fin}}_w(H_A(\beta))$
is called a {\it fake exponent}.

Let $w\cdot u$ denote
$w_1u_1+\cdots +w_nu_n$ for $u\in {\bf Q}^n$.
Let 
\begin{equation}
\label{eqn:L}
L=\{\, u\in {\bf Z}^n\, :\, Au=0\,\}.
\end{equation}
An $A$-hypergeometric series
\begin{equation}
\label{eq:SeriesInW}
x^v\cdot\sum_{u\in L}
g_u(\log x)x^u
\qquad (g_u\in {\bf k}[x])
\end{equation}
is said to be {\it in the direction of $w$} if
there exists
a basis $u^{(1)},\ldots, u^{(n)}$ of
${\bf Q}^n$ with $w\cdot u^{(j)}>0$ 
($j=1,\ldots,n$) such that
$g_u=0$ whenever $u\notin \sum_{j=1}^n {\bf Q}_{\geq 0}u^{(j)}$.
A fake exponent $v$ is called an {\it exponent}
if there exists
an $A$-hypergeometric series (\ref{eq:SeriesInW}) 
in the direction of $w$ 
with nonzero $g_0$.
Let $\prec$ be the lexicographic order on ${\bf N}^n$.
Suppose that $u^{(1)},\ldots, u^{(n)}$ is a basis as above.
Then a monomial like $x^v\cdot {\rm in}_\prec (g_0)(\log x)$ 
in the $A$-hypergeometric series 
\begin{equation}
\label{eq:CanonicalSeries}
x^v\cdot\sum_{u\in L\cap \sum_{j=1}^n {\bf Q}_{\geq 0}u^{(j)}}
g_u(\log x)x^u
\qquad (g_u\in {\bf k}[x])
\end{equation}
with nonzero $g_0$
is called a
{\it starting monomial}.
The $A$-hypergeometric 
series (\ref{eq:CanonicalSeries})
is said to be {\it canonical} with respect to
$w$ if no starting monomials other than
$x^v\cdot {\rm in}_\prec (g_0)(\log x)$ appear in the series.

Next we recall logarithm-free $A$-hypergeometric series $\phi_v$.
For $v\in {\bf k}^n$, its {\it negative support} 
${\rm nsupp}(v)$ is the set of indices
$i$ with $v_i\in {\bf Z}_{<0}$.
When ${\rm nsupp}(v)$ is minimal with respect to inclusions
among ${\rm nsupp}(v+u)$ with $u\in L$,
$v$ is said to have {\it minimal negative support}.
For $v$ satisfying $Av=\beta$ with minimal negative support,
we define a formal series
\begin{equation}
\label{eqn:CanonicalSeries}
\phi_v=x^v \psi_v=x^v\cdot
\sum_{u\in N_v}\frac{[v]_{u_-}}{[v+u]_{u_+}}x^{u}.
\end{equation}
Here 
\begin{equation}
N_v=\{\, u\in L\, :\,
 {\rm nsupp}(v)={\rm nsupp}(v+u)\,\},
\end{equation}
and
$u_+, u_-\in {\bf N}^n$ satisfy $u=u_+ -u_-$ with disjoint supports,
and
$[v]_t=\prod_{j=1}^n v_j(v_j-1)\cdots (v_j-t_j+1)$
for $t\in {\bf N}^n$.
Proposition 3.4.13 and Theorem 3.4.14 in \cite{sst-book}
respectively state that
the series $\phi_v$ is $A$-hypergeometric, and that
if $v$ is a fake exponent of $M_A(\beta)$, then
$\phi_v$ is canonical, and $v$ is an exponent.
Let 
${\rm Minex}_{\beta,w}$ denote the set of fake exponents with
minimal negative support, of $M_A(\beta)$ with respect to $w$.
Then ${\rm Minex}_{\beta,w}$ is a subset of the set of exponents;
moreover by Corollary 3.4.15 in \cite{sst-book},
under the correspondence of $v$ with $\phi_v$,
the set ${\rm Minex}_{\beta,w}$ corresponds to the set of
logarithm-free canonical series solutions to
$M_A(\beta)$ with respect to $w$.
Let ${\cal S}_{\beta,w}$ denote the space spanned by logarithm-free $A$-hypergeometric
series in the direction of $w$.
Then the set of logarithm-free canonical $A$-hypergeometric series
with respect to $w$ is a basis of ${\cal S}_{\beta,w}$.
Hence we have
\begin{equation}
\label{eq:Sw}
\dim {\cal S}_{\beta,w}=|{\rm Minex}_{\beta,w}|.
\end{equation}
We shall compute $|{\rm Minex}_{\beta,w}|$ in the following sections.

\section{Finite sets $E_\tau(\beta)$ and 
logarithm-free canonical $A$-hypergeometric series}

We denote the set $\{\, a_1,\ldots, a_n\,\}$ of the column vectors
of the matrix $A$ by $A$ as well.
Let ${\rm conv}(A)$ denote the convex hull of the set $A$, and
$\Delta_w$ the regular
triangulation of ${\rm conv}(A)$ determined by $w$.
In this section, we define a finite set $E_\tau(\beta)$
associated to a parameter $\beta$ and a face $\tau$ of $\Delta_w$,
and prove that each element of $E_\tau(\beta)$ gives
volume-of-$\tau$-many logarithm-free canonical series solutions
with respect to $w$.

Let ${\bf N}A$ denote the monoid generated by $A$.
Given a face $\tau$ of $\Delta_w$, ${\bf Z}(A\cap\tau)$
and ${\bf k}(A\cap\tau)$ respectively denote
the additive group and the vector space generated by 
$A\cap\tau$.
Here
we agree that ${\bf k}(A\cap\tau)={\bf Z}(A\cap\tau)=\{\, 0\,\}$
when $\tau=\emptyset$.
Associated to a parameter $\beta\in {\bf k}^d$
and a face $\tau\in \Delta_w$, we define a finite set
$E_\tau(\beta)$ in the same way as in \cite{IsoClass}:
\begin{equation}
\label{eqn:EtauBeta}
E_\tau (\beta):=
\{\, \lambda \in {\bf k}(A\cap\tau) /{\bf Z}(A\cap\tau)\, :\,
\beta-\lambda \in {\bf N}A+{\bf Z}(A\cap\tau)\,\}.
\end{equation}
An algorithm for $E_\tau(\beta)$ is given in \cite{Saito-Traves}.
For fundamental properties  of $E_\tau(\beta)$ such as the finiteness,
see \cite{IsoClass}.

Let ${\rm vert}(\tau)$ denote the set of vertices
of a face $\tau$.
For a face $\tau$ of $\Delta_w$ or
for a subset $\tau$ of $A$, we abuse the notation to
denote the set $\{\, j\, :\, a_j\in\tau\,\}$ by $\tau$ again;
for example, the symbol
${\rm vert}(\tau)$ could mean
the set $\{\, j\, :\, \mbox{$a_j$ is a vertex of $\tau$}\,\}$.
For a face $\tau$ of $\Delta_w$, let 
${\rm vol}(\tau)$ denote the index 
$[{\bf Z}(A\cap \tau) : \sum_{i\in {\rm vert}(\tau)}{\bf Z}a_i]$.

\begin{theorem}
\label{theorem:GeneralConvergence}
Let $w$ be a generic weight vector.
Given a face $\tau\in \Delta_w$,
there exist at least $|E_\tau(\beta)| {\rm vol}(\tau)$-many
logarithm-free canonical series solutions $\phi_v$
to $M_A(\beta)$ with respect to $w$,
satisfying $\sum_{j\in\tau}v_ja_j\in E_\tau(\beta)$,
${\rm nsupp}(v)\subset {\rm vert}(\tau)$, and $\{\, j\, :\, v_j\notin {\bf Z}\,\}
\subset {\rm vert}(\tau)$.
\end{theorem}

{\bf Proof.}
Put $I:={\rm vert}(\tau)$.
Given $\lambda\in E_\tau(\beta)$,
there exists $v\in {\bf k}^n$ such that
$\beta=Av$, $\lambda =\sum_{a_j\in\tau}v_ja_j$,
${\rm nsupp}(v)\subset \tau$, and
$\{\, j\, :\, v_j\notin {\bf Z}\,\}\subset \tau$.
Since $I$ is a base of $\tau$,
we may further assume ${\rm nsupp}(v)\subset I$
and $\{\, j\, :\, v_j\notin {\bf Z}\,\}\subset I$.
There are exactly ${\rm vol}(\tau)$-many such vectors $v$ modulo
$L$ (see (\ref{eqn:L}) for the definition of $L$),
which we prove below as Lemma \ref{lemma:IndexNumber}.
We can take the representatives $v$ such that
they have minimal negative supports.
Then the series $\phi_v$ 
(see (\ref{eqn:CanonicalSeries}) for the definition of $\phi_v$) 
are formal solutions, and
$$
N_v=\{\, u\in L\, :\, 
u_j<-v_j\, (j\in {\rm nsupp}(v)),\,\,
u_j\geq -v_j\, (j\in {\rm psupp}(v))\,\}.
$$
Here, similarly to the negative support,
${\rm psupp}(v)$ denotes
the set of indices $i$ with $v_i\in {\bf N}$.

Let $L_{{\bf Q}}$ be the ${\bf Q}$-vector space generated by $L$,
and consider a polyhedron
$$
N_{v,{\bf Q}}:=\{\, u\in L_{{\bf Q}}\, :\, 
u_j\leq -v_j\, (j\in {\rm nsupp}(v)),\,\,
u_j\geq -v_j\, (j\in {\rm psupp}(v))\,\}.
$$
Since the vectors $a_i$ ($i\in I$) are linearly independent,
the characteristic cone of $N_{v,{\bf Q}}$
$$
K_{v,I}=\{\, u\in L_{{\bf Q}}\, :\, 
u_j\leq 0\, (j\in {\rm nsupp}(v)),\,\,
u_j\geq 0\, (j\in {\rm psupp}(v))\,\}
$$
is pointed, or strongly convex.
Hence the polyhedron $N_{v,{\bf Q}}$ has vertices and
it is decomposed as 
$
N_{v,{\bf Q}}= P +K_{v,I},
$
where $P$ is the convex hull of the set of vertices of $N_{v,{\bf Q}}$
(for example, see \cite{Schrijver}).
Note that
the characteristic cone $K_{v, I}$ of $N_{v,{\bf Q}}$ for any such $v$
is contained in the cone $K_I:=
\{ \, u\in L_{\bf Q}\, :\, u_j\geq 0\quad (j\notin I)\,\}.$
Let $\tilde{\tau}\in\Delta_w$ be a $d$-dimensional face with
$\tilde{\tau}\supset\tau$.
Then $I\subset {\rm vert}(\tilde{\tau})=:\tilde{I}$, and
\begin{eqnarray}
\label{eqn:CharacteristicCone}
K_I\subset K_{\tilde{I}}
&=&\{ \, u\in L_{\bf Q}\, :\, u_j\geq 0\quad (j\notin \tilde{I})\,\}
\nonumber\\
&=& \sum_{j\notin \tilde{I}}
{\bf Q}_{\geq 0}(e_j-\sum_{i\in \tilde{I}}c_{ji}e_i),
\end{eqnarray}
where $c_{ji}$ are given by $a_j=\sum_{i\in \tilde{I}}c_{ji}a_i$,
and $e_1,\ldots,e_n$ is the standard basis for ${\bf Z}^n$.
Since $\tilde{\tau}\in \Delta_w$ means that $w\cdot u>0$ for all nonzero 
$u\in K_{\tilde{I}}$,
and since $w$ is generic,
there exists a unique optimal solution $u$ 
for the integer programming problem:
$$
\mbox{Minimize $w\cdot u$ subject to $u\in N_v$.}
$$
Replacing $v$ by $v+u$,
we may assume $0$ is the unique optimal solution.
Then we claim that $v$ is a fake exponent, i.e.,
$x^v$ is annihilated by $\sum_{j=1}^n a_{ij}\theta_j -\beta_i$
($i=1,\ldots, d$) and ${\rm in}_w(I_A)$.
Since $Av=\beta$, we only need to verify that $x^v$ is annihilated by
${\rm in}_w(I_A)$.
Let $Au_+=Au_-$ and $w\cdot u_+ > w\cdot u_-$.
Among all the terms appearing in 
$\partial^{u_+}(\phi_v)-\partial^{u_-}(\phi_v)$,
$\partial^{u_+}(x^v)=[v]_{u_+}x^{v-u_+}$ is the unique term 
with cost $w\cdot (v-u_+)$.
Since $\phi_v$ is a formal solution, it implies $\partial^{u_+}(x^v)=0$.
We have thus proved that $v$ is a fake exponent.
It is furthermore an exponent, and $\phi_v$ is canonical
by Theorem 3.4.14 in \cite{sst-book}.
\qed

For a subset $I$ of $\{\, 1,\ldots, n\,\}$, let
$I^c$ denote the complement of $I$, i.e., $I^c=\{\, 1,\ldots, n\,\}
\setminus I$.
Given $\lambda\in E_\tau(\beta)$, put
\begin{eqnarray}
{\cal E}xp_w(\tau,\lambda)
:=
\{\, v\in {\bf k}^n &:&
\beta = Av,\,\, \lambda=\sum_{j\in\tau}v_ja_j,\nonumber\\
&&v_j\in {\bf N}\,\,\mbox{for all $j\in {\rm vert}(\tau)^c$}\,\}.
\end{eqnarray}
We denote by ${\cal E}xp_w(\tau,\lambda)/L$
the image of ${\cal E}xp_w(\tau,\lambda)$
in $\{\, v\in {\bf k}^n\, :\, \beta=Av\,\}/L$.

\begin{lemma}
\label{lemma:IndexNumber}
We keep the notation in the proof of 
Theorem \ref{theorem:GeneralConvergence}.
Let $\lambda\in E_\tau(\beta)$. Then
\begin{eqnarray}
\label{eqn:eqn0}
&&{\bf Z}(A\cap\tau)/\sum_{i\in I}{\bf Z}a_i\\
\label{eqn:eqn1}
&\simeq&  
\{\, u\in L_{\tau, {\bf Q}}\, :
\, u_i\in {\bf Z}\,\,\mbox{for all $i\in \tau\setminus I$}\,\}/L_\tau
\\
&\simeq &
\label{eqn:eqn2}
{\cal E}xp_w(\tau,\lambda)/L,
\end{eqnarray}
where $
L_\tau:=\{\, u\in L\, :\,
u_j=0\,\,\mbox{for all $j\in\tau^c$}\,\}$, and
$L_{\tau, {\bf Q}}$ is the ${\bf Q}$-vector space generated by $L_\tau$.
\end{lemma}

{\bf Proof.}
As we saw in the proof of Theorem \ref{theorem:GeneralConvergence},
the set (\ref{eqn:eqn2}) is not empty.

First we prove the bijectivity between the sets 
(\ref{eqn:eqn0}) and (\ref{eqn:eqn1}).
For $u\in L_{\tau, {\bf Q}}$ with $u_i\in {\bf Z}$ 
($i\in \tau\setminus I$), we define 
$g(u):=\sum_{i\in \tau\setminus I}u_ia_i\in {\bf Z}(A\cap\tau)$.
For the surjectivity, suppose $u_i\in {\bf Z}$ for all $i\in \tau\setminus I$.
Since $I$ is a base of $\tau$, we can find $b_i\in {\bf Q}$ ($i\in I$)
such that  $\sum_{i\in \tau\setminus I}u_ia_i = \sum_{i\in I}b_ia_i$.
Then $\sum_{i\in \tau\setminus I}u_ia_i=g(\tilde{u})$,
where $\tilde{u}_i=u_i$ for $i\in \tau\setminus I$ and
$\tilde{u}_i=-b_i$ for $i\in I$.
For the injectivity of $g$,
suppose $g(u)=g(v)$.
Then $\sum_{i\in I}u_ia_i-\sum_{i\in I}v_ia_i=
\sum_{i\in \tau\setminus I}v_ia_i-\sum_{i\in \tau\setminus I}u_ia_i
\in \sum_{i\in I}{\bf Z}a_i$.
Since $a_i$ ($i\in I$) are linearly independent,
$u_i-v_i$ is an integer for all $i\in I$.
Hence $u-v\in L$.

Next let $u$ belong to the set (\ref{eqn:eqn1}).
Clearly we may assume that $u_i\in {\bf N}$ for all $i\in \tau\setminus I$.
Suppose that $v$ is an element of the set (\ref{eqn:eqn2}).
Then $v+u$ is again an element of the set (\ref{eqn:eqn2}).
Thus the set (\ref{eqn:eqn1})
can be embedded into the set (\ref{eqn:eqn2}).

Next we show that the set (\ref{eqn:eqn2})
can be embedded into the set (\ref{eqn:eqn0}).
Suppose that $v$ and $v'$ belong to the set (\ref{eqn:eqn2}).
Then $\sum_{i\in I}(v_i -v'_i)a_i\in {\bf Z}(A\cap\tau)$.
If it belongs to $\sum_{i\in I}{\bf Z}a_i$, then $v_i-v'_i\in {\bf Z}$
for all $i\in I$. Hence $v-v'\in L$.

Finally we remark that the set (\ref{eqn:eqn0})
is finite, to complete the proof.
\qed

\section{Dimension formula for the vector space of 
logarithm-free $A$-hypergeometric
series}

In this section, we give a dimension formula for the vector space ${\cal S}_{\beta,w}$.
Recall that
the set ${\rm Minex}_{\beta,w}$
corresponds to the set of
logarithm-free canonical series solutions to
$M_A(\beta)$ with respect to $w$, and hence
$\dim {\cal S}_{\beta,w} = |{\rm Minex}_{\beta,w}|$.

We denote by ${\rm vert}(\Delta_w)$
the set of vertices of $\Delta_w$.
\begin{lemma}
\label{lemma:FakeExponent1}
If $v$ is a fake exponent, then
$v_j\in {\bf N}$ for $j\notin {\rm vert}(\Delta_w)$.
\end{lemma}

{\bf Proof.}
If $a_j\notin {\rm vert}(\Delta_w)$, then
there exist nontrivial $m_j, m_{ji}\in {\bf N}$ such that
$m_ja_j =\sum_{i\in {\rm vert}(\Delta_w)}m_{ji}a_i$.
Hence we have $\partial_j^{m_j}\in {\rm in}_w(I_A)$
for $j\notin {\rm vert}(\Delta_w)$.
Since $v$ is a fake exponent, we see $v_j\in {\bf N}$ for $j\notin 
{\rm vert}(\Delta_w)$.
\qed

Suppose $v\in {\rm Minex}_{\beta,w}$. 
Define a set $I_v$ by
\begin{equation}
I_v:=\{\, j\, :\,
v_j\notin {\bf N}\,\}.
\end{equation}
By Lemma \ref{lemma:FakeExponent1}, 
the set $I_v$ is a subset of ${\rm vert}(\Delta_w)$.
Since $v$ has minimal negative support,
$I_v$ is linearly independent, i.e.,
$a_j$ ($j\in I_v$) are linearly independent.

\begin{lemma}
\label{lemma:Lambda}
Let $v\in {\rm Minex}_{\beta,w}$, $I_v:=\{\, j\, :\,
v_j\notin {\bf N}\,\}$, and
$\tau_v:={\rm conv}(\{ a_j\, :\, j\in I_v \})$.
Then $\tau_v$
is a face of $\Delta_w$,
$I_v$ is the vertex set of $\tau_v$, and
$\sum_{a_i\in \tau_v}
v_ia_i$ represents an element of $E_{\tau_v}(\beta)$.
\end{lemma}

{\bf Proof.}
We prove only the statement $\tau_v\in \Delta_w$.
The others follow easily.
Suppose $\tau_v\notin \Delta_w$.
Then there exist a subset $I\subset I_v$,
integers $m_i\in {\bf N}$ ($i\in I$) and
$n_j\in {\bf N}$ ($j\notin I$) with the properties
$\sum_{i\in I}m_i a_i =\sum_{j\notin I} n_j a_j$ and
$\sum_{i\in I}m_i w_i > \sum_{j\notin I} n_j w_j$.
Since $v$ is a fake exponent,
$\prod_{i\in I}\partial_i^{m_i}\in {\rm in}_w(I_A)$
implies $v_i\in {\bf N}$ for some $i\in I$.
This contradicts the definition of $I_v$.
\qed

Lemma \ref{lemma:Lambda} shows that $v\in {\rm Minex}_{\beta,w}$ implies
$v\in {\cal E}xp_w(\tau_v,\lambda_v)$, where
\begin{equation}
\tau_v:={\rm conv}(\{ a_j\, :\, j\in I_v \})\quad\mbox{and}\quad
\lambda_v:=\sum_{a_i\in \tau_v}
v_ia_i\in E_{\tau_v}(\beta).
\end{equation}
Define a map $\Lambda$ from ${\rm Minex}_{\beta,w}$ to the set
$
\coprod_{\tau\in\Delta_w, \lambda\in E_\tau(\beta)}
{\cal E}xp_w(\tau,\lambda)/L
$
by 
\begin{equation}
\Lambda(v)=[v]\in {\cal E}xp_w(\tau_v,\lambda_v)/L,
\end{equation}
where $[v]$ denotes the image of $v$ 
in $\{\, v\in {\bf k}^n\, :\, \beta=Av\,\}/L$.

\begin{lemma}
\label{lemma:Fiber}
Let $v,v'\in {\rm Minex}_{\beta,w}$. Suppose $\tau_v=\tau_{v'}$ and
$v+L=v'+L$.
Then $v=v'$.
\end{lemma}

{\bf Proof.}
Since $I_v$ (resp. $I_{v'}$) is the vertex set of 
$\tau_v$ (resp. $\tau_{v'}$),
we have $I_v=I_{v'}$, which implies ${\rm psupp}(v)={\rm psupp}(v')$.
Then by the equality $v+L=v'+L$, we have
${\rm nsupp}(v)={\rm nsupp}(v')$.
If $v\not= v'$, then
without loss of generality, we may assume $w\cdot (v -v')<0$.
Since $\phi_v$ is canonical by Theorem 3.4.14 in \cite{sst-book},
the term $x^{v'}$ does not appear in $\phi_v$, which contradicts the 
definition of $\phi_v$.
\qed

The following is immediate from Lemma \ref{lemma:Fiber}.

\begin{corollary}
\label{cor:injectivity}
The map $\Lambda$ is injective, and
\begin{equation}
|{\rm Minex}_{\beta,w}|=\sum_{\tau\in\Delta_w, \lambda\in E_\tau(\beta)}
|\Lambda(\Lambda^{-1}(
{\cal E}xp_w(\tau,\lambda)/L))|.
\end{equation}
\end{corollary}

We introduce a partial order $\preceq$ in the set
$
\{ \, (\tau,\lambda)\, :\, \tau\in\Delta_w,\, \lambda\in E_\tau(\beta)\,\}$.
Note that there exists a natural map from $E_{\tau'}(\beta)$
to $E_\tau(\beta)$ if $\tau'\preceq\tau$, i.e., if $\tau'$ is a face of $\tau$.
We define a partial order $\preceq$ by 
$(\tau',\lambda')\preceq (\tau,\lambda)$
if $\tau'\preceq\tau$, and
the image of $\lambda'$ under the natural map from $E_{\tau'}(\beta)$
to $E_\tau(\beta)$ coincides with $\lambda$.
Note that ${\cal E}xp_w(\tau',\lambda')/L$ is a subset of
${\cal E}xp_w(\tau,\lambda)/L$ 
when $(\tau',\lambda')\preceq (\tau,\lambda)$.

\begin{lemma}
\label{ImageOfLambda}
\begin{equation}
\Lambda(\Lambda^{-1}(
{\cal E}xp_w(\tau,\lambda)/L))
={\cal E}xp_w(\tau,\lambda)/L\setminus
\bigcup_{(\tau',\lambda')\prec (\tau,\lambda)}
{\cal E}xp_w(\tau',\lambda')/L.
\end{equation}
\end{lemma}

{\bf Proof.}
The definition of $\Lambda$ implies that the LHS is contained in
the RHS.
Suppose that $c$ belongs to the RHS.
By the proof of Theorem \ref{theorem:GeneralConvergence}, there exists
$v\in {\rm Minex}_{\beta,w}\cap {\cal E}xp_w(\tau,\lambda)$ 
such that $c=[v]$ in $\{\, v\, :\, \beta=Av\,\}/L$.
Since $c$ does not belong to
$\bigcup_{(\tau',\lambda')\prec (\tau,\lambda)}
{\cal E}xp_w(\tau',\lambda')/L$,
neither does $[v]$, which means $\Lambda(v)\in {\cal E}xp_w(\tau,\lambda)$,
and thus $\Lambda(v)=c$.
\qed

For a face $\tau$ of $\Delta_w$, let ${\rm facet}(\tau)$ denote
the set of one-codimensional faces of $\tau$.
The following is the main theorem in this section.

\begin{theorem}
\label{theorem:LogFreeFormula}
\begin{enumerate}
\item
Let $\tau\in\Delta_w$ and $\lambda\in E_\tau(\beta)$. Then
\begin{eqnarray}
&&|\Lambda(\Lambda^{-1}(
{\cal E}xp_w(\tau,\lambda)/L))|
\nonumber\\
&& =
{\rm vol}(\tau) -
\sum_{i}\sum_{\stackrel{\scriptstyle E_{\tau_i}(\beta)\ni\lambda}
{\tau_1,\ldots,\tau_i\in {\rm facet}(\tau)}}
(-1)^{i-1}{\rm vol}(\tau_1\cap\cdots\cap\tau_i).
\end{eqnarray}
\item
\begin{eqnarray}
\dim {\cal S}_{\beta,w} &=&
\sum_{\tau\in\Delta_w} \left( \sum_{\lambda\in E_\tau(\beta)}
{\rm vol}(\tau)\right.
\nonumber\\
&-&\left.
\sum_{i}\sum_{\stackrel{\scriptstyle E_{\tau_i}(\beta)\ni\lambda}
{\tau_1,\ldots,\tau_i\in {\rm facet}(\tau)}}
(-1)^{i-1}{\rm vol}(\tau_1\cap\cdots\cap\tau_i)
\right).
\end{eqnarray}
\end{enumerate}
\end{theorem}

{\bf Proof.}
It is enough to prove the first statement by Corollary \ref{cor:injectivity}.
The combination of
Lemma \ref{lemma:IndexNumber}, Corollary \ref{cor:injectivity}, and
Lemma \ref{ImageOfLambda} implies that
the number $|\Lambda^{-1}({\cal E}xp_w(\tau,\lambda)/L)|$ equals
the number of the equivalence classes of 
${\bf Z}(A\cap\tau)/\sum_{i\in {\rm vert}(\tau)}{\bf Z}a_i$
which cannot be represented by an element of
those for $(\tau',\lambda')\prec (\tau,\lambda)$.
An equivalence class of ${\bf Z}(A\cap\tau)/\sum_{i\in 
{\rm vert}(\tau)}{\bf Z}a_i$ which can be represented by an element of
${\bf Z}(A\cap\tau')$ for $(\tau', \lambda')\prec (\tau,\lambda)$,
is actually represented by one for $\tau'$ of codimension one to $\tau$.
The number of such equivalence classes for $\tau'$ is
${\rm vol}(\tau')$.
Here we remark that ${\bf Z}(A\cap \tau')\cap {\bf Z} ({\rm vert}(\tau))
= {\bf Z} ({\rm vert}(\tau'))$ since ${\rm vert}(\tau)$ is linearly independent.
An equivalence class which can be represented by equivalence classes for
$\tau_1,\ldots, \tau_i$ comes from an equivalence class for
$\tau_1\cap\cdots\cap\tau_i$.
Hence the first statement holds.
\qed

\begin{example}
\label{example:SmallExample}
Let
$$ A \quad = \quad
\pmatrix{ 
1 & 1 & 1 \cr
0 & 1 & 2 }, $$ 
and $\beta=a_1={}^t(1,1)$.
Then there are two regular triangulations $\Delta_w$ and
$\Delta_{w'}$ of ${\rm conv}(A)$
(see Figure \ref{figure:triangulation}), and
for all $\tau$ of either triangulation,
$E_\tau(\beta)=\{ 0\}$.
\begin{figure}[h]
\begin{picture}(40,60)(-130,10)
\put(20,10){\circle*{5}}
\put(20,30){\circle*{5}}
\put(20,10){\line(0,1){20}}
\put(10,10){\vector(1,0){30}}
\put(10,10){\vector(1,2){20}}
\put(40,30){$\Delta_w$}
\put(110,10){\circle*{5}}
\put(110,20){\circle*{5}}
\put(110,30){\circle*{5}}
\put(110,10){\line(0,1){20}}
\put(100,10){\vector(1,0){30}}
\put(100,10){\vector(1,2){20}}
\put(130,30){$\Delta_{w'}$}
\end{picture}
\caption{$\Delta_w$ and $\Delta_{w'}$}
\label{figure:triangulation}
\end{figure}
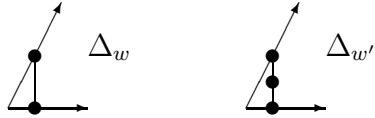

By Theorem \ref{theorem:LogFreeFormula}, we obtain
\begin{eqnarray}
\dim {\cal S}_{\beta,w} &=& (2-2\times 1 +1)+2(1-1) +1
\nonumber\\
&=& 1+0+0+1=2,
\end{eqnarray}
and
\begin{eqnarray}
\dim {\cal S}_{\beta, w'} &=& 2(1-2\times 1 +1)+3(1-1) +1
\nonumber\\
&=& 0+0+0+1=1.
\end{eqnarray}
In both cases, the polynomial solution $\phi_{(0,1,0)}=x_2$
is associated with face $\emptyset$.
In the case $\Delta_w$, there exists a logarithm-free 
$A$-hypergeometric series $\phi_{(1/2,0,1/2)}$
associated with face ${\rm conv}(A)$.
\end{example}

A generalization of the case $\Delta_w$ in Example
\ref{example:SmallExample} will be considered in Section
\ref{section:SimplicialCase}.
The case $\Delta_{w'}$ can be generalized as the
following corollary.
A triangulation $\Delta_w$ is said to be {\it unimodular}
if $a_j$ ($j\in {\rm vert}(\tau)$)
is a basis of ${\bf Z}A$ for all $d$-dimensional faces $\tau\in\Delta_w$.

\begin{corollary}
Suppose that $\beta\in {\bf N}A$, and
$\Delta_w$ is unimodular.
Then ${\cal S}_{\beta,w}$ is one-dimensional,
and coincides with the space of polynomial solutions.
\end{corollary}

{\bf Proof.}
In this case, ${\rm vol}(\tau)=1$ and $E_\tau(\beta)=\{ 0\}$ for all
faces $\tau\in \Delta_w$.
For a face $\tau\not=\emptyset$, we have
$$
{\rm vol}(\tau) -
\sum_{i}\sum_{\tau_1,\ldots,\tau_i\in {\rm facet}(\tau)}
(-1)^{i-1}{\rm vol}(\tau_1\cap\cdots\cap\tau_i)=0.
$$
The face $\emptyset$ gives one canonical series.
It is a polynomial.
\qed

\section{Logarithm
coefficients of $A$-hypergeometric series}

In the next section, we shall prove that
there exists a fundamental system of solutions to $M_A(\beta)$
consisting of logarithm-free canonical series with respect to
a suitable weight $w$ in the case where ${\rm conv}(A)$ is a simplex.
In order to prove this, we consider logarithm
coefficients of $A$-hypergeometric series
in this section.

We denote by ${\bf k}[L]$
the symmetric algebra of the vector space
spanned by
$\sum_{j=1}^n c_jx_j$ ($c\in L$),
which is a subring of the polynomial ring 
${\bf k}[x] = {\bf k}[x_1,\ldots, x_n]$.

\begin{lemma}
\label{lemma:VanishingOfA}
Let $f\in {\bf k}[x]$.
Then
$\sum_{j=1}^n a_{ij}\partial_j (f) =0$ for all $i=1,\ldots, d$
if and only if
$f\in {\bf k}[L]$.
\end{lemma}

{\bf Proof.}
We prove the only-if-direction.
The if-direction is straightforward.
Without loss of generality,
we may assume that the vectors $a_1,\ldots, a_d$
are linearly independent.
Then ${\bf k}[x]={\bf k}[L][x_1,\ldots, x_d]$.
Let $f=\sum_v g_v x^v$ with $g_v\in {\bf k}[L]$
and $v\in {\bf N}^d$.
Suppose that there exists $u\not= 0$ with
$g_u\not= 0$.
Without loss of generality we may assume that $u_1\not=0$.
Then by the equality 
$\sum_{j=1}^n a_{ij}\partial_j (f) =0$
we have
$\sum_{j=1}^d a_{ij}(u_j+1-\delta_{1j})g_{u-e_1+e_j}=0$
for all $i$,
where $e_1,\ldots,e_d$ is the standard basis for ${\bf Z}^d$.
This contradicts the linear independence of the vectors
$a_1,\ldots, a_d$.
\qed

\begin{proposition}
\label{proposition:LogCoef}
Suppose that 
\begin{equation}
\phi=\sum_{u} g_u( \log (x))x^u
\end{equation}
is an $A$-hypergeometric series.
Then $g_u\in {\bf k}[L]$ for all $u$.
\end{proposition}

{\bf Proof.}
Suppose that $\phi$ is a solution to $M_A(\beta)$.
Then
$(\sum_{j=1}^na_{ij}\theta_j-\beta_i)(g_u( \log (x))x^u)=0$,
which leads to $Au=\beta$ and
$\sum_{j=1}^n a_{ij}\partial_j (g_u) =0$.
Lemma \ref{lemma:VanishingOfA} finishes the proof.
\qed

Let $v\in {\bf k}$ and $n\in {\bf N}$.
For $j\in {\bf N}$ with $j\leq n$,
put
$$
c_j(v,n):=\sum_{1\leq i_1<\cdots<i_j\leq n}
\prod_{k=1}^j(v-i_k+1).
$$
For example, $c_n(v,n)=v(v-1)\cdots (v-n+1)$,
$c_{n-1}(v,n)=\sum_{i=1}^n\prod_{k\not= i}(v-k+1)$,
$c_1(v,n)=\sum_{i=1}^n (v-i+1) =nv -n(n-1)/2$.
We agree that $c_0(v,n)=1$.
Then $c_{n-j}(v,n)=(1/j!)(d/dv)^j(c_n(v,n))$.
Note that for $n>0$
\begin{equation}
\label{eqn:c_{n-1}}
c_{n-1}(v,n)\not= 0\qquad
\mbox{if $v\in {\bf N}$.}
\end{equation}

Let $v\in {\bf k}^n$ and $u\in {\bf N}^n$.
For $u'\in {\bf N}^n$ with $u'\leq u$, i.e., with
$u'_j\leq u_j$ for all $j$,
put
$$
c_{u'}(v,u):=\prod_{j=1}^n c_{u'_j}(v_j,u_j).
$$
The proof of the following lemma is straightforward.

\begin{lemma}
\label{lemma:LogCoef2}
Let $p\in {\bf k}[x]$.
Then
$$
\partial^u(x^v p(\log x))
=\sum_{0\leq u'\leq u}
c_{u'}(v,u)x^{v-u} (\partial^{u-u'}p)(\log x).
$$
\end{lemma}

\begin{proposition}
\label{proposition:LogCoef2}
Let $w$ be a generic weight vector.
Let $v$ be an exponent of $M_A(\beta)$ 
with respect to $w$.
Assume that an $A$-hypergeometric series
$$
\phi=\sum_{u\in L}g_u(\log x)x^{v+u}
$$
satisfies $g_u=0$ for $w\cdot u<0$.
Then $g_u(x)\in
{\bf k}[\, x_j\, :\, j\in {\rm vert}(\Delta_w)\,]$
for all $u\in L$.
\end{proposition}

{\bf Proof.}
Put $V:={\rm vert}(\Delta_w)$, and
let $j\notin V$.
As in the proof of Lemma \ref{lemma:FakeExponent1},
there exist $m_j, m_{ji}\in {\bf N}$ with the properties
$m_ja_j =\sum_{i\in V}m_{ji}a_i$
and $w_jm_j> \sum_{i\in V}w_im_{ji}$.
Suppose that there exists $u\in L$ such that
the degree $d_j$ of $g_u$ in $x_j$ is positive.
Let $N$ be a positive integer, and $h$ a polynomial
satisfying
$$
\partial_j^{Nm_j}(g_u(\log x)x^{v+u})
=h(\log x) x^{v+u-Nm_j e_j}.
$$
By (\ref{eqn:c_{n-1}}) and Lemma \ref{lemma:LogCoef2},
we see that the degree of $h$ in $x_j$ is at least $d_j-1$.
Hence
the coefficient of $x^{v+u-Nm_j e_j+\sum_{i\in V}Nm_{ji}e_i}$
in $\phi$ is not zero.
For large $N$, 
$w\cdot (u-Nm_j e_j+\sum_{i\in V}Nm_{ji}e_i)<0$.
This contradicts our assumption.
\qed

\section{Simplicial case}
\label{section:SimplicialCase}

Throughout this section, we assume that the convex hull
${\rm conv}(A)$ is a simplex,
and that $w$ is a generic weight vector 
such that the corresponding regular
triangulation $\Delta_w$ is ${\rm conv}(A)$ itself.
Then $\tau\mapsto {\bf Q}_{\geq 0}\tau$ gives a 
one-to-one correspondence between the faces of $\Delta_w$ and those of
the cone
\begin{equation}
{\bf Q}_{\geq 0}A=
\{\, \sum_{j=1}^n c_ja_j\, :\, c_j\in {\bf Q}_{\geq 0}\,\}.
\end{equation}
Hence we consider the same sets $E_\tau(\beta)$ as the ones in \cite{IsoClass}.
In this section, first we prove that all canonical series solutions with
respect to $w$ are logarithm-free.
Second we give a rank formula
for $M_A(\beta)$ in terms of the finite sets $E_\tau(\beta)$.
Third we characterize the exceptional set 
in terms of $E_\tau(\beta)$ again.
Finally we prove the equivalence of the Cohen-Macaulayness of the affine
toric variety defined by $A$ with the emptiness of the exceptional set.

Let $V$ denote the set of vertices of ${\rm conv}(A)$, i.e.,
$V={\rm vert}({\rm conv}(A))$.
Recall that a fake exponent is a zero
of the fake indicial ideal (see (\ref{def:fin}) for the definition).

\begin{proposition}
\label{proposition:FakeIndicial}
The fake indicial ideal 
$\widetilde{{\rm fin}}_w(H_A(\beta))$
is radical.
\end{proposition}

{\bf Proof.}
Let $v$ be a fake exponent.
It is enough to prove that
the localization $\widetilde{{\rm fin}}_w(H_A(\beta))_v$
at $v$ is radical.
As seen in the proof of Lemma \ref{lemma:FakeExponent1},
if $j\notin V$, then
there exists $m_j\in {\bf N}$ such that
$\partial_j^{m_j}\in {\rm in}_w(I_A)$.
Hence $\theta_j-v_j\in \widetilde{{\rm fin}}_w(H_A(\beta))_v$
for $j\notin V$.
Since each component of $A\theta-\beta$ belongs to
$\widetilde{{\rm fin}}_w(H_A(\beta))$,
we see that each component of
$\sum_{j\in V}(\theta_j -v_j)a_j$ belongs to
$\widetilde{{\rm fin}}_w(H_A(\beta))_v$.
Then we conclude $\theta_j -v_j\in \widetilde{{\rm fin}}_w(H_A(\beta))_v$
for $j\in V$ as well
because the square matrix $( a_j: j\in V)$ is nonsingular.
\qed

\begin{proposition}
\label{proposition:Exponents}
A fake exponent which does not have minimal negative
support is never an exponent.
\end{proposition}

{\bf Proof.}
Let $v$ be an exponent and $v'$ a fake exponent
such that $v-v'\in L$, and that
${\rm nsupp}(v)$ is a proper subset of ${\rm nsupp}(v')$.
Put $u:=v-v'$.
We first prove that $w\cdot u_{+}> w\cdot u_{-}$.
Assume the contrary.
Since $v'$ is a fake exponent,
$\partial^{u_-}(x^{v'})=0$.
Hence there exists $i$ such that
$v'_i\in {\bf N}$ and
$v'_i< -u_i=v'_i-v_i$.
This implies $i\in {\rm nsupp}(v)$, and thus
$i\in {\rm nsupp}(v')$.
This contradicts the fact that $v'_i\in {\bf N}$.

Suppose that $\phi$ is a canonical solution with
starting monomial $x^{v'}$.
Since $V={\rm vert}({\rm conv}(A))$
is linearly independent, $\phi$ is logarithm-free by
Propositions \ref{proposition:LogCoef} and \ref{proposition:LogCoef2}.
As seen in the first paragraph, $\partial^{u_-}(x^{v'})\not=0$.
Hence the coefficient of $x^v$ in $\phi$ is not zero.
This contradicts the fact that $\phi$ is canonical.
\qed

By Propositions \ref{proposition:FakeIndicial}, 
and \ref{proposition:Exponents},
Theorem \ref{theorem:LogFreeFormula} gives a rank formula in this case.

\begin{theorem}
\label{simplex:log-free-sol}
Assume that ${\rm conv}(A)$ is a simplex.
Let $w$ be a generic weight vector such that 
the corresponding regular triangulation $\Delta_w$ is 
${\rm conv}(A)$ itself.
Then there exists a fundametal system of solutions to $M_A(\beta)$
consisting of logarithm-free canonical series
with respect to $w$, and the rank of $M_A(\beta)$, denoted by 
${\rm rank}(M_A(\beta))$, is given by the following formula:
\begin{eqnarray}
{\rm rank}(M_A(\beta)) &=&
\sum_{\tau} \left( \sum_{\lambda\in E_\tau(\beta)}
{\rm vol}(\tau)\right.
\nonumber\\
\label{equation:Formula}
&-&\left.
\sum_{i}\sum_{\stackrel{\scriptstyle E_{\tau_i}(\beta)\ni\lambda}
{\tau_1,\ldots,\tau_i\in {\rm facet}(\tau)}}
(-1)^{i-1}{\rm vol}(\tau_1\cap\cdots\cap\tau_i)
\right).
\end{eqnarray}
\end{theorem}

\begin{remark}
Theorem \ref{simplex:log-free-sol} generalizes
Theorem 4.2.4 and Corollary 4.2.5 in \cite{sst-book},
which are concerned with the case $d=2$.
\end{remark}

For a generic parameter $\beta$, the rank of the $A$-hypergeometric system
$M_A(\beta)$ is known to equal the volume ${\rm vol}(A)$
(see \cite{Adolphson} and \cite{sst-book}).
We characterize such a generic parameter by 
considering the following equivalent conditions on $\beta$:

\begin{condition}
\label{condition:Generic}
\begin{enumerate}
\item
If $\lambda\in {\bf k}(A\cap\tau_1 \cap \cdots\cap \tau_l)$
belongs to $E_{\tau_i}(\beta)$ for all $i=1,\ldots, l$, then
$\lambda\in E_{\tau_1\cap\cdots\cap\tau_l}(\beta)$.
\item
If $\lambda\in {\bf k}(A\cap\tau_1 \cap \tau_2)$
belongs to $E_{\tau_i}(\beta)$ for $i=1,2$, then
$\lambda\in E_{\tau_1\cap\tau_2}(\beta)$.
\end{enumerate}
\end{condition}

The set
$$
{\cal E}(A):=
\{\, \beta\, :\, 
{\rm rank}(M_A(\beta))>{\rm vol}(A)
\,\}
$$
is called the {\it exceptional set}.

\begin{theorem}
\label{theorem:ExceptionalSet}
The exceptional set ${\cal E}(A)$
is the set of parameters which
does not satisfy 
Condition \ref{condition:Generic}.
\end{theorem}

{\bf Proof.}
Note that ${\rm rank}(M_A(\beta))=|{\rm Minex}_{\beta,w}|$
by Propositions \ref{proposition:FakeIndicial} and
\ref{proposition:Exponents}, and ${\rm vol}(A)=
|{\bf Z}A/\sum_{j\in V}{\bf Z}a_j|$.
Recall the proof of Theorem \ref{theorem:LogFreeFormula}.
Given an exponent $v\in {\rm Minex}_{\beta,w}$, we put
$I_v:=\{\, j\, :\, v_j\notin {\bf N}\,\}\subset V$,
$\tau_v:=\sum_{j\in I_v}{\bf Q}_{\geq 0} a_j$,
and
$\lambda_v:=\sum_{j\in I_v}v_ja_j\in E_{\tau_v}(\beta)$.
Fix an exponent $v(0)\in {\rm Minex}_{\beta,w}$, and
define a map 
$C:{\rm Minex}_{\beta,w}\rightarrow {\bf Z}A/\sum_{j\in V}{\bf Z}a_j$
by 
$$
C(v)=\sum_{j\in V}(v_j-v(0)_j)a_j =\sum_{j\notin V}(v_j-v(0)_j)a_j.
$$
Then $C$ is surjective by Theorem \ref{theorem:GeneralConvergence}.

We show that Condition \ref{condition:Generic} is equivalent to
the injectivity of $C$.
Suppose that $\beta$ satisfies
Condition \ref{condition:Generic}, and
let $v(1),v(2)\in {\rm Minex}_{\beta,w}$
satisfy $C(v(1))=C(v(2))$.
Then $\sum_{j\in V}(v(1)_j-v(2)_j)a_j\in \sum_{j\in V}{\bf Z}a_j$.
The linear independence of $V$ implies that
$v(1)-v(2)\in L$.
Put $\lambda:=\sum_{j\in I_{v(1)}\cap I_{v(2)}}v(1)_ja_j\in 
{\bf k}(A\cap\tau_{v(1)}\cap\tau_{v(2)})$.
Then $\lambda =\lambda_{v(i)}$ in $E_{\tau_{v(i)}}(\beta)$ ($i=1,2$),
and thus 
by
Condition \ref{condition:Generic},
$\lambda\in E_{\tau_{v(1)}\cap\tau_{v(2)}}(\beta)$.
By Theorem \ref{theorem:GeneralConvergence},
there exists $v\in {\rm Minex}_{\beta, w}$
such that $I_v\subset I_{v(1)}\cap I_{v(2)}$,
and $\sum_{j\in I_{v(1)}\cap I_{v(2)}}(v_j-v(1)_j)a_j
\in \sum_{j\in I_{v(1)}\cap I_{v(2)}}{\bf Z}a_j$.
Then $C(v)=C(v(1))=C(v(2))$.
This leads to the fact that $v-v(1), v-v(2)\in L$ as above,
and thus ${\rm nsupp}(v)\subset {\rm nsupp}(v(i))$ ($i=1,2$).
Since $v(1)$, $v(2)$, and $v$ have minimal negative support,
we obtain $v(1)=v=v(2)$.

Next we assume the injectivity of $C$, and suppose that 
$\lambda\in {\bf k}(A\cap\tau_1\cap\tau_2)$ belongs
to $E_{\tau_1}(\beta)$ and $E_{\tau_2}(\beta)$.
By Theorem \ref{theorem:GeneralConvergence},
there exists $v(i)\in {\rm Minex}_{\beta, w}$ such that
$I_{v(i)}\subset {\rm vert}(\tau_i)$, and
$\lambda -\sum_{j\in I_{v(i)}}v(i)_j a_j 
\in \sum_{j\in {\rm vert}(\tau_i)}{\bf Z}a_j$
($i=1, 2$).
Hence $C(v(1))=C((v(2))$, and thus $v(1)=v(2)$
by the injectivity.
Put $v:=v(1)=v(2)$.
Then we have $I_v\subset {\rm vert}(\tau_1)\cap {\rm vert}(\tau_2)$,
and
$\beta-\lambda= \sum_{j\notin I_v} v_ja_j +
(\sum_{j\in I_v} v_ja_j-\lambda)\in {\bf N}A + 
\sum_{j\in {\rm vert}(\tau_1)\cap {\rm vert}(\tau_2)}{\bf Z}a_j$.
Hence we obtain $\lambda\in E_{\tau_1\cap\tau_2}(\beta)$.
\qed

\begin{theorem}
\label{theorem:CMequiv}
The ring ${\bf k}[{\bf N}A]={\bf k}[\partial]/I_A$ is Cohen-Macaulay 
if and only if
${\cal E}(A)=\emptyset$.
\end{theorem}

{\bf Proof.}
Since $\{\, \partial_j\,\}_{j\in V}$ is a linear system of
parameters, the Cohen-Macaulayness means that ${\bf k}[{\bf N}A]$
is a free ${\bf k}[\partial_j\, :\, j\in V]$-module.
Hence ${\bf k}[{\bf N}A]$ is not Cohen-Macaulay
if and only if there exist
minimal generators $\partial^{u(1)}$ and $\partial^{u(2)}$,
and $\partial^{m(1)}, \partial^{m(2)}\in {\bf k}[\partial_j\, :\,
j\in V]$ with ${\rm supp}(m(1))\cap {\rm supp}(m(2))=\emptyset$
such that $\partial^{m(1)}\partial^{u(1)}-
\partial^{m(2)}\partial^{u(2)}\in I_A$, where ${\rm supp}(m(i))$
denotes the support of $m(i)$.
Suppose that ${\bf k}[{\bf N}A]={\bf k}[\partial]/I_A$ is not
Cohen-Macaulay.
Let 
$\tau_i:=\sum_{j\in {\rm supp}(m(i))}{\bf Q}_{\geq 0}a_j$
for $i=1,2$, and
$\beta:=Au(1)-Am(2)=Au(2)-Am(1)$.
Then
$0\in E_{\tau_i}(\beta)$ ($i=1,2$)
and $0\notin E_{\{\, 0\,\}}(\beta)$,
and hence $\beta\in {\cal E}(A)$.

Replacing $\beta$ by $\beta-\lambda$ 
in Condition \ref{condition:Generic} if necessary,
we see that
the condition ${\cal E}(A)\not=\emptyset$ is equivalent to the condition
that there exist $\beta\in {\bf Z}A$, and faces $\tau_1,\tau_2$
 such that $0\in E_{\tau_i}(\beta)$
($i=1,2$) and $0\notin E_{\tau_1\cap\tau_2}(\beta)$.
The latter condition is equivalent to the condition 
that there exist $\beta\in {\bf Z}A$, and
$m(1), m(2)\in {\bf N}^V$
such that $\beta+A m(i)\in {\bf N}A$ ($i=1,2$),
and no $m(3)\in {\bf N}^V$ with
${\rm supp}(m(3))\subset {\rm supp}(m(1))\cap{\rm supp}(m(2))$
such that $\beta+A m(3)\in {\bf N}A$.
Let $m'$ be the greatest common divisor of $m(1)$ and $m(2)$.
By considering $\beta +Am'$ instead of $\beta$,
we see that the above condition is equivalent
to the condition that
there exist $\beta\in {\bf Z}A\setminus {\bf N}A$, and
$m(1), m(2)\in {\bf N}^V$
with ${\rm supp}(m(1))\cap{\rm supp}(m(2))=\emptyset$
such that $\beta+A m(i)\in {\bf N}A$ ($i=1,2$).

Suppose ${\cal E}(A)\not= \emptyset$.
Let $\beta\in {\bf Z}A\setminus {\bf N}A$, and
$m(1), m(2)\in {\bf N}^V$ be the ones in the last condition.
Then for $i=1,2$, there exist $u(i)\in {\bf N}^n$ and $v(i)\in {\bf N}^V$ 
such that $\partial^{u(i)}$ is a minimal generator,
and $\beta +Am(i)=Av(i)+Au(i)$.
We may assume that ${\rm supp}(m(i))\cap{\rm supp}(v(i))=\emptyset$
($i=1,2$).
Then $\partial^{m(2)+v(1)}\partial^{u(1)}-
\partial^{m(1)+v(2)}\partial^{u(2)}\in I_A$.
Since the supports of $m(2)+v(1)$ and $m(1)+v(2)$ are different,
$\partial^{u(1)}\not= \partial^{u(2)}$.
We conclude that
${\bf k}[{\bf N}A]={\bf k}[\partial]/I_A$ is not
Cohen-Macaulay.
\qed

\begin{remark}
In the general situation, 
the only-if-direction of Theorem \ref{theorem:CMequiv} 
was proved by
Gel'fand, Zelevinskii, and Kapranov (\cite{GZK},
see also \cite{Adolphson}).
The if-direction was conjectured by Sturmfels, and proved when $d=2$
(see \cite{CDD}).
Recently Matusevich \cite{Laura} proved it when $n-d=2$.
\end{remark}

\begin{example}
Let
$$ A \quad = \quad
\pmatrix{ 
1 & 1 & 1 & 1 & 1 & 1 & 1 & 1 & 1 \cr
0 & 1 & 2 & 3 & 0 & 2 & 0 & 1 & 0 \cr
0 & 0 & 0 & 0 & 1 & 1 & 2 & 2 & 3}, $$ 
and
$\beta_0= {}^t(1,1,1)$.
We prove that ${\cal E}(A) = \{\, \beta_0\,\}$ and
${\rm rank}(M_A(\beta_0))= 11$.
Note that ${\bf N}A$ has clear ${\bf Z}/3{\bf Z}$-symmetry 
(see Figure \ref{figure:rank11}).

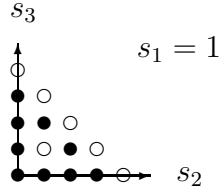
\begin{figure}[ht]
\begin{picture}(50,60)(-140,10)
\put(10,10){\circle*{5}}
\put(20,10){\circle*{5}}
\put(30,10){\circle*{5}}
\put(40,10){\circle*{5}}
\put(50,10){\circle{5}}
\put(10,20){\circle*{5}}
\put(20,20){\circle{5}}
\put(30,20){\circle*{5}}
\put(40,20){\circle{5}}
\put(10,30){\circle*{5}}
\put(20,30){\circle*{5}}
\put(30,30){\circle{5}}
\put(10,40){\circle*{5}}
\put(20,40){\circle{5}}
\put(10,50){\circle{5}}
\put(10,10){\vector(1,0){50}}
\put(10,10){\vector(0,1){50}}
\put(7,70){$s_3$}
\put(70,8){$s_2$}
\put(55,55){$s_1=1$}
\end{picture}
\caption{The set $A$}
\label{figure:rank11}
\end{figure}

Given a subset $\{\, i_1,\ldots, i_r\}$ of $\{\, 1,4,9\,\}$,
let $\tau_{i_1\cdots i_r}$ denote the face 
${\bf Q}_{\geq 0}a_{i_1}+\cdots +{\bf Q}_{\geq 0}a_{i_r}$ of
${\bf Q}_{\geq 0}A$.
Then there are one $3$-dimensional face $\tau_{149}={\bf Q}_{\geq 0}A$,
three $2$-dimesional faces $\tau_{14}, \tau_{19}, \tau_{49}$,
three $1$-dimensional faces $\tau_{1}, \tau_{4}, \tau_{9}$,
and one $0$-dimensional face $\tau_\emptyset=\{\, 0\,\}$.
The volumes corresponding to the faces are 
\begin{eqnarray}
&&{\rm vol}(\tau_{149})={\rm vol}(A)=9,\nonumber\\
&&{\rm vol}(\tau_{14})={\rm vol}(\tau_{19})={\rm vol}(\tau_{49})
=3,\nonumber\\
&&{\rm vol}(\tau_{1})={\rm vol}(\tau_{4})={\rm vol}(\tau_{9})
=1,\nonumber\\
&&{\rm vol}(\{\, 0\,\})=1.
\end{eqnarray}

\noindent
For each face $\tau$, we have 
$({\bf Q}(A\cap\tau))\cap {\bf Z}A={\bf Z}(A\cap\tau)$,
and thus $|E_\tau(\beta)|=0$ or $1$.
Given $\beta$, let ${\rm minface}(\beta)$ denote the set of
minimal faces among faces $\tau$ with nonempty $E_\tau(\beta)$.
Then the equivalence class of $M_A(\beta)$ is determined by
${\rm minface}(\beta)$ (see Theorem 2.1 in \cite{IsoClass}).

Suppose $\beta\in {\cal E}(A)$.
Theorem \ref{theorem:ExceptionalSet} 
implies that the set ${\rm minface}(\beta)$
contains at least two elements.
If ${\rm minface}(\beta)$ contains $\tau_1,\tau_4$,
then it has to contain $\tau_9$ as well.
Hence ${\rm minface}(\beta)=\{\, \tau_1,\tau_4,\tau_9\,\}$,
$\beta=\beta_0$, and
\begin{eqnarray}
{\rm rank}(M_A(\beta_0))&=& (9-3\times 3 +3\times 1 -1)
\nonumber\\
&+& 3(3-2\times 1 +1)
\nonumber\\
&+& 3\times 1= 2 + 6 + 3 =11 > 9={\rm vol}(A).
\end{eqnarray}
If $E_{\tau_{14}}(\beta)\not=\emptyset$ and 
$E_{\tau_{19}}(\beta)\not=\emptyset$,
then we see $E_{\tau_{1}}(\beta)\not=\emptyset$.
Hence the cases ${\rm minface}(\beta)=\{\, \tau_{14},\tau_{19}\,\}$,
$\{\, \tau_{14},\tau_{19},\tau_{49}\,\}$, and
$\{\, \tau_{1},\tau_{49}\,\}$ do not occur.
This finishes the proof of the equality ${\cal E}(A) = \{\, \beta_0\,\}$.

\end{example}

\section{Analytic Classification}

We gave a criterion for two $A$-hypergeometric systems to be isomorphic
as algebraic $D$-modules in \cite{IsoClass}.
In this section, by using the argument in the proof of 
Theorem \ref{theorem:GeneralConvergence},
we show that the same criterion remains valid in the analytic
category.

\begin{lemma}
\label{2EtauBeta}
Let $w$ be a generic weight vector. 
Let $\tau$ be a face of ${\bf Q}_{\geq 0}A$, and
$I$ the vertex set ${\rm vert}(\tau')$ of a face $\tau'\in \Delta_w$
with
${\bf Q}\tau = \sum_{i\in I}{\bf Q}a_i$.
Then there exists a natural map 
\begin{equation}
\nu: E_{\tau'}(\beta)\longrightarrow E_\tau(\beta).
\end{equation}
The map $\nu$ is surjective and a
$[{\bf Z}(A\cap\tau) : {\bf Z}(A\cap\tau')] : 1$-map,
unless both $E_{\tau'}(\beta)$ and $E_\tau(\beta)$ are empty.
In particular, 
\begin{equation}
|E_{\tau'}(\beta)|=
|E_{\tau}(\beta)|[{\bf Z}(A\cap\tau) : {\bf Z}(A\cap\tau')].
\end{equation}
\end{lemma}

{\bf Proof.}
First we show the existence of $\nu$.
Note that ${\bf k}(A\cap\tau)={\bf k}(A\cap\tau')$, and
let $\lambda$ belong to ${\bf k}(A\cap\tau)={\bf k}(A\cap\tau')$.
Then $\lambda\in E_{\tau'}(\beta)$ means
$\beta-\lambda\in {\bf N}A +{\bf Z}(A\cap\tau')$,
which implies 
$\beta-\lambda\in {\bf N}A +{\bf Z}(A\cap\tau)$,
which in turn means $\lambda\in E_\tau(\beta)$.
This gives the natural map $\nu$.

Second we prove the surjectivity of $\nu$.
Suppose $\lambda\in E_\tau(\beta)$.
Then as in the first two sentences of the proof of Theorem
\ref{theorem:GeneralConvergence},
there exists $v\in {\bf k}^n$ such that
$\beta=Av$, $\lambda =\sum_{a_j\in\tau}v_ja_j$,
${\rm nsupp}(v)\subset \tau'$, and
$\{\, j\, :\, v_j\notin {\bf Z}\,\}\subset \tau'$.
Then $\sum_{a_j\in\tau'}v_ja_j\in E_{\tau'}(\beta)$
and $\nu(\sum_{a_j\in\tau'}v_ja_j)=\lambda$.

Finally we prove that $\nu$ is a 
$[{\bf Z}(A\cap\tau) : {\bf Z}(A\cap\tau')] : 1$-map.
Suppose that $\lambda,\lambda'\in E_{\tau'}(\beta)$ and 
$\nu(\lambda)=\nu(\lambda')$.
Then $\lambda-\lambda'\in {\bf Z}(A\cap\tau)$.
Conversely suppose $\lambda\in E_{\tau'}(\beta)$, and
$v\in {\bf Z}^n$ with $v_j=0$ for $j\notin\tau$.
We prove that $\lambda-Av\in E_{\tau'}(\beta)$.
Since $I$ is a base of $\tau$,
there exists $l\in L_\tau$ such that
$(v+l)_j\in {\bf N}$  for all $j\notin I$.
Then $\beta-(\lambda - A(v+l))\in {\bf N}A +{\bf Z}(A\cap\tau')$.
Hence $\lambda-Av=\lambda -A(v+l)\in E_{\tau'}(\beta)$.
\qed

\begin{theorem}
\label{theorem:convergence}
Let $w$ be a generic weight vector.
Let $\tau$ be a face of ${\bf Q}_{\geq 0}A$, and
$I$ the vertex set ${\rm vert}(\tau')$ of a face $\tau'\in \Delta_w$
with
${\bf Q}\tau = \sum_{i\in I}{\bf Q}a_i$.
Then
there exist at least
$|E_\tau(\beta)| {\rm vol}_{\tau}(\tau')$-many
linearly independent logarithm-free canonical series solutions $\phi_v$
to $M_A(\beta)$ with respect to $w$,
satisfying $\sum_{a_j\in\tau}v_ja_j\in E_\tau(\beta)$,
${\rm nsupp}(v)\subset I$, and $\{\, j\, :\, v_j\notin {\bf Z}\,\}
\subset I$.
Here ${\rm vol}_{\tau}(\tau')$ denotes the index
$[{\bf Z}(A\cap \tau) : \sum_{i\in I}{\bf Z}a_i]$.\end{theorem}

{\bf Proof.}
This is immediate from Theorem \ref{theorem:GeneralConvergence}
and Lemma \ref{2EtauBeta}.
\qed

The algebraic version of the following is the main theorem 
in the paper \cite{IsoClass}.

\begin{theorem}
\label{thm:AnalyticVersion}
Let ${\bf k}$ be the field of complex numbers ${\bf C}$,
and ${\cal D}$ the sheaf of rings of
analytic linear differential operators on ${\bf C}^n$.
Put ${\cal M}_A(\beta)={\cal D}\otimes_D M_A(\beta)$,
and ${\cal M}_A(\beta')$ likewise. 
The $A$-hypergeometric systems ${\cal M}_A(\beta)$ and ${\cal M}_A(\beta')$ are
isomorphic as ${\cal D}$-modules if and only if
$E_\tau(\beta)=E_\tau(\beta')$ for all faces $\tau$ of the cone
${\bf Q}_{\geq 0}A$.
\end{theorem}

{\bf Proof.}
Here we prove the only-if-part of the theorem.
The if-part is immediate from Theorem 2.1 in \cite{IsoClass}.

We suppose
that $\lambda\in E_\tau(\beta)\setminus E_\tau(\beta')$
for some face $\tau$,
and then we prove that ${\cal M}_A(\beta)$ and 
${\cal M}_A(\beta')$ are not isomorphic.

Take $w$, $v$ with $\lambda=\sum_{j\in\tau}v_j a_j$,
and $I$ as in Theorem 
\ref{theorem:convergence}.
For $j\in\tau\setminus I$, let $u(j)=e_j-\sum_{i\in I}c_{ji}e_i$
in the notation of (\ref{eqn:CharacteristicCone}).
Since $\tau$ is a face of ${\bf Q}_{\geq 0}A$, and since
$I$ is a base of $\tau$,
we have $K_I=\sum_{j\in\tau\setminus I}{\bf Q}_{\geq 0} u(j)$,
and
$u(j)$ ($j\in\tau\setminus I$) is a basis of $L_{\tau,{\bf Q}}$
satisfying $L\cap K_I \subset \sum_{j\in \tau\setminus I}{\bf N}u(j)$.
The series $\psi_{v}$ (see (\ref{eqn:CanonicalSeries}))
converges in $U_r$ for $r$ small enough, where
$$
U_r:=\{\, x\, :\, |x^{u(j)}|< r \,\, (j\in\tau\setminus I),
\quad x_i\not= 0\,\, (i\in \tau)\,\},
$$
since $M_A(\beta)$ is regular holonomic
(for example, see Theorem 2.4.9 in \cite{sst-book}).
Hence $\phi_v$ is a nonzero element of
${\rm Hom}_{{\cal D}(U_r)}({\cal M}_A(\beta)(U_r), x^v{\cal O}(U_r))$.
Here note that $x^v {\cal O}_{U_r}$ is a ${\cal D}_{U_r}$-module.

Next we prove that
${\rm Hom}_{{\cal D}(U_r)}({\cal M}_A(\beta')(U_r), x^v{\cal O}(U_r))=0$,
which completes the proof of the theorem.
Suppose that $\phi=x^v\psi$ belongs to the solution space
${\rm Hom}_{{\cal D}(U_r)}({\cal M}_A(\beta')(U_r), x^v{\cal O}(U_r))$.
Since $\psi$ is holomorphic on $U_r$,
it has the Laurent expasion
\begin{equation}
\psi (x) = \sum_{l\in {\bf N}^{\tau^c}\times {\bf Z}^\tau}
c_l y^l,
\end{equation}
where $y_j=x^{u(j)}$ if $j\in \tau\setminus I$, and $y_j=x_j$ otherwise.
By the equality $(A\theta -\beta')x^v\psi =0$,
each $l\in {\bf N}^{\tau^c}\times {\bf Z}^\tau$ with $c_l\not=0$
must satisfy $\beta'=\sum_{i\in I}(v_i+l_i)a_i + 
\sum_{j\in \tau\setminus I} v_ja_j+
\sum_{j\notin \tau} (v_j+l_j)a_j$.
The condition that $\lambda\notin E_\tau(\beta')$,
implies that all $c_l$ are zero.
Therefore ${\cal M}_A(\beta)$ and ${\cal M}_A(\beta')$ are not isomorphic.
\qed


\bigskip

Department of Mathematics

Hokkaido University

Sapporo, 060-0810

Japan

e-mail: saito@math.sci.hokudai.ac.jp

\end{document}